\documentclass{article}
\usepackage{amssymb}
\usepackage{amsmath}
\usepackage{chicago}

\newtheorem{theorem}{Theorem}

\newtheorem{corollary}[theorem]{Corollary}

\newtheorem{lemma}[theorem]{Lemma}

\newtheorem{proposition}[theorem]{Proposition}

\begin{document}

\title{Berry-Esseen for Free Random Variables}
\author{Vladislav Kargin \thanks{%
Courant Institute of Mathematical Sciences; 109-20 71st Road, Apt. 4A, New
York NY 11375; kargin@cims.nyu.edu} \thanks{%
I am grateful to Diana Bloom for her help with editing this paper. }}
\date{}
\maketitle

\begin{abstract}
An analogue of the Berry-Esseen inequality is proved for the speed of
convergence of free additive convolutions of bounded probability measures.
The obtained rate of convergence is of the order $n^{-1/2},$ the same as in
the classical case. An example with \ binomial measures shows that this
estimate cannot be improved without imposing further restrictions on
convolved measures.
\end{abstract}

\section{Introduction}

In the setting of non-commutative probability theory, an interesting new
operation on probability measures was defined. If two probability measures, $%
\mu_{1}$ and $\mu_{2},$ are represented as spectral measures of two free
self-adjoint operators, $X_{1}$ and $X_{2},$ then the sum $X_{1}+X_{2}$ has
the spectral probability measure $\mu_{3}$, which is called\emph{\ free
convolution} of $\mu_{1}$ and $\mu_{2},$ and denoted as $\mu_{1}\boxplus
\mu_{2}$. The theory of additive free convolutions was invented by %
\shortciteN{voiculescu83}, \citeyear{voiculescu86}, and an excellent account
of further development of free probability theory can be found in %
\shortciteN{hiai_petz00}. Free convolution has many properties similar to
the usual convolution of probability measures. In particular, an analogue of
the Central Limit Theorem holds, which says that multiple free convolutions
of a probability measure with itself converge to the Wigner semicircle law.
In this paper we investigate the speed of convergence and establish an
inequality similar to the classical Berry-Esseen inequality.

Without reference to operator theory, the free convolution can be defined as
follows. Suppose $\mu_{1}$ and $\mu_{2}$ are two probability measures
compactly supported on the real line. Define the \emph{Cauchy transform} of $%
\mu_{i}$ as 
\begin{equation*}
G_{i}\left( z\right) =\int_{-\infty}^{\infty}\frac{d\mu_{i}\left( t\right) }{%
z-t}. 
\end{equation*}
Each of $G_{i}\left( z\right) $ is well-defined and univalent for large
enough $z$ and we can define its functional inverse, which is well-defined
in a neighborhood of 0. Let us call this inverse the $K$\emph{-function} of $%
\mu_{i}$ and denote it as $K_{i}\left( z\right) $:%
\begin{equation*}
K_{i}\left( G_{i}\left( z\right) \right) =G_{i}\left( K_{i}\left( z\right)
\right) =z. 
\end{equation*}
Then, we define $K_{3}\left( z\right) $ by the following formula:%
\begin{equation}
K_{3}\left( z\right) =K_{1}\left( z\right) +K_{2}\left( z\right) -\frac{1}{z}%
.  \label{summation_property}
\end{equation}
It turns out that $K_{3}\left( z\right) $ is the $K$-function of a
probability measure, $\mu_{3}$, which is the\emph{\ free convolution} of $%
\mu_{1}$ and $\mu_{2}$.

Let us turn to issues of convergence. Let $d\left( \mu_{1},\mu_{2}\right) $
denote the Kolmogorov distance between the probability measures $\mu_{1}$
and $\mu_{2}.$ That is, if $\mathcal{F}_{1}$ and $\mathcal{F}_{2}$ are the
distribution functions corresponding to measures $\mu_{1}$ and $\mu_{2}$,
respectively, then 
\begin{equation*}
d\left( \mu_{1},\mu_{2}\right) =:\sup_{x\in\mathbb{R}}\left\vert \mathcal{F}%
_{1}\left( x\right) -\mathcal{F}_{2}\left( x\right) \right\vert . 
\end{equation*}
Let $\mu$ be a probability measure with the zero mean and unit variance and
let $m_{3}$ be its third absolute moment. Then the classical Berry-Esseen
inequality says that 
\begin{equation*}
d\left( \mu^{\left( n\right) },\nu\right) \leq Cm_{3}\frac{1}{\sqrt{n}}, 
\end{equation*}
where $\nu$ is the standard Gaussian measure and $\mu^{\left( n\right) }$ is
the normalized $n$-time convolution of measure $\mu$ with itself:%
\begin{equation*}
\mu^{\left( n\right) }\left( du\right) =\mu\ast...\ast\mu\left( \sqrt {n}%
du\right) . 
\end{equation*}
This inequality was proved by \shortciteN{berry41} and \shortciteN{esseen45}
for a more general situation of independent but not necessarily identical
measures. A simple example with binomial measures shows that in this
inequality the order of $n^{-1/2}$ cannot be improved without further
restrictions.

We aim to derive a similar inequality when the usual convolution of measures
is replaced by free convolution. Namely, let 
\begin{equation*}
\mu^{\left( n\right) }\left( du\right) =\mu\boxplus...\boxplus\mu\left( 
\sqrt{n}du\right) 
\end{equation*}
and let $\nu$ denote the standard semicircle distribution. It is known that $%
\mu^{\left( n\right) }$ converges weakly to $\nu$ (\shortciteN{voiculescu83}%
, \shortciteN{maassen92}, \shortciteN{pata96}, and \shortciteN{voiculescu98}%
). We are interested in the speed of this convergence and we prove that if $%
\mu$ is supported on $\left[ -L,L\right] , $ then%
\begin{equation}
d\left( \mu^{\left( n\right) },\nu\right) \leq CL^{3}\frac{1}{\sqrt{n}}.
\label{inequality_main}
\end{equation}
An example shows that the rate of $n^{-1/2}$ cannot be improved without
further restrictions$,$ similar to the classical case.

The main tool in our proof of inequality (\ref{inequality_main}) is Bai's
theorem (\citeyearNP{bai93}) that relates the supremum distance between two
probability measures to a distance between their Cauchy transforms. To
estimate the distance between Cauchy transforms, we use the fact that as $n$
grows, the $K$-function of $\mu^{\left( n\right) }$ approaches the $K$%
-function of the semicircle law. Therefore, the main problem in our case is
to investigate whether the small distance between $K$-functions implies a
small distance between the Cauchy transforms themselves. We approach this
problem using the Lagrange formula for functional inverses.

The rest of the paper is organized as follows: Section 2 contains the
formulation and the proof of the main result. It consists of several
subsections. In Subsection 2.1 we formulate the result and outline its
proof. Subsection 2.2 evaluates how fast the $K$-function of $\mu^{\left(
n\right) }$ approaches the $K$-function of the semicircle law. Subsection
2.3 provides useful estimates on behavior of the Cauchy transform of the
semicircle law and related functions. Subsection 2.4 introduces a functional
equation for the Cauchy transforms and concludes the proof by estimating how
fast the Cauchy transform of $\mu^{\left( n\right) }$ converges to the
Cauchy transform of the semicircle law. An example in Section 3 shows that
the rate of $n^{-1/2}$ cannot be improved.

\section{Formulation and Proof of the Main Result}

\subsection{Formulation and outline of the proof}

Let the \emph{semicircle law} be the probability measure on the real line
that has the following cumulative distribution function:%
\begin{equation*}
\Phi\left( x\right) =\frac{1}{\pi}\int_{-\infty}^{x}\sqrt{4-t^{2}}\chi_{%
\left[ -2,2\right] }\left( t\right) \,dt, 
\end{equation*}
where $\chi_{\left[ -2,2\right] }\left( t\right) $ is the characteristic
function of the interval $\left[ -2,2\right] .$

\begin{theorem}
\label{theorem_main}Suppose that $\mu$ is a probability measure that has
zero mean and unit variance, and is supported on interval $\left[ -L,L\right]
$. Let $\mu^{\left( n\right) }$ be the normalized $n$-time free convolution
of measure $\mu$ with itself: $\mu^{\left( n\right) }\left( du\right)
=\mu\boxplus...\boxplus\mu\left( \sqrt{n}du\right) .$ Let $\mathcal{F}%
_{n}\left( x\right) $ denote the cumulative distribution function of $%
\mu^{\left( n\right) }$. Then for large enough $n$ the following bound
holds: 
\begin{equation*}
\sup_{x}\left\vert \mathcal{F}_{n}\left( x\right) -\Phi\left( x\right)
\right\vert \leq CL^{3}n^{-1/2}, 
\end{equation*}
where $C$ is an absolute constant.
\end{theorem}

\textbf{Remark:} $C=2^{16}$ will do, although this constant is far from the
best possible.

\textbf{Proof: }\ First, we quote one of Bai's results:

\begin{theorem}[Bai (1993)]
Let measures with distribution functions $\mathcal{F}$ and $\mathcal{G}$ be
supported on a finite interval $\left[ -B,B\right] $ and let their Cauchy
transforms be $G_{\mathcal{F}}\left( z\right) $ and $G_{\mathcal{G}}\left(
z\right) ,$ respectively. Let $z=u+iv.$ Then 
\begin{align}
\sup_{x}\left\vert \mathcal{F}\left( x\right) -\mathcal{G}\left( x\right)
\right\vert & \leq\frac{1}{\pi\left( 1-\kappa\right) \left( 2\gamma
-1\right) }\left[ \int_{-A}^{A}\left\vert G_{\mathcal{F}}\left( z\right) -G_{%
\mathcal{G}}\left( z\right) \right\vert du\right.  \label{Bai_inequality} \\
& \left. +\frac{1}{v}\sup_{x}\int_{\left\vert y\right\vert
\leq2vc}\left\vert \mathcal{G}\left( x+y\right) -\mathcal{G}(x)\right\vert dy%
\right] ,  \notag
\end{align}
where $A>B,$ $\kappa=\frac{4B}{\pi\left( A-B\right) \left( 2\gamma -1\right) 
}<1,$ $\gamma>1/2,$ and $c$ and $\gamma$ are related by the following
equality 
\begin{equation*}
\gamma=\frac{1}{\pi}\int_{\left\vert u\right\vert <c}\frac{1}{1+u^{2}}du. 
\end{equation*}
\end{theorem}

Note that if $\mathcal{G}\left( x\right) $ is the semicircle distribution
then $\left\vert \mathcal{G}^{\prime}\left( x\right) \right\vert \leq
\pi^{-1}.$ Therefore $\left\vert \mathcal{G}\left( x+y\right) -\mathcal{G}%
(x)\right\vert \leq\left\vert y\right\vert /\pi.$ $\ $Integrating this
inequality, we obtain:%
\begin{equation}
\frac{1}{v}\sup_{x}\int_{\left\vert y\right\vert \leq2vc}\left\vert \mathcal{%
G}\left( x+y\right) -\mathcal{G}(x)\right\vert dy\leq\frac{4c^{2}}{\pi}v.
\label{inequality_semicircle_smoothness}
\end{equation}
Hence, the main question is how fast $v$ can be made to approach zero if the
first integral in (\ref{Bai_inequality}) is also required to approach zero.

Let $G_{\Phi}$ and $G_{n}$ be the Cauchy transforms of the semicircle law
and $\mu^{\left( n\right) },$ respectively. Assume for the moment that the
following lemma holds:

\begin{lemma}
\label{Cauchy_tr_closeness copy(1)}Suppose $v=2^{10}L^{3}/\sqrt{n}.$ Then
for all sufficiently large $n,$ we have the following estimate:%
\begin{equation*}
\int\limits_{-8}^{8}\left\vert G_{n}\left( u+iv\right) -G_{\Phi}\left(
u+iv\right) \right\vert du\leq\frac{2^{10}L^{3}}{\sqrt{n}}. 
\end{equation*}
\end{lemma}

Let us apply Bai's theorem using Lemma \ref{Cauchy_tr_closeness copy(1)} and
inequality (\ref{inequality_semicircle_smoothness}). Let $\Phi\left(
x\right) $ and $\mathcal{F}_{n}\left( x\right) $ denote the cumulative
distribution functions of the semicircle law and $\mu^{\left( n\right) },$
respectively. The semicircle law is supported on $\left[ -2,2\right] ,$ and
for any fixed interval $I$ that includes $\left[ -2,2\right] ,$ we can find
such $n_{0}$ that $\mu^{\left( n\right) }$ is supported on $I$ for all $%
n\geq n_{0}$ (see \shortciteN{bercovici_voiculescu95}). Suppose that $n$ is
so large that $\mu^{\left( n\right) }$ is supported on $\left[
-2^{5/4},2^{5/4}\right] .$ Then we can take $A=8,$ $B=2^{5/4},$ $c=6$, and
calculate $\gamma=0.895$ and $\kappa=0.682$. Then Bai's theorem gives the
following estimate: 
\begin{align*}
\sup_{x}\left\vert \mathcal{F}_{n}\left( x\right) -\Phi\left( x\right)
\right\vert & \leq1.268\left( 2^{10}L^{3}n^{-1/2}+46937L^{3}n^{-1/2}\right)
\\
& \leq2^{16}L^{3}n^{-1/2}.
\end{align*}
QED.

Thus, the main task is to prove Lemma \ref{Cauchy_tr_closeness copy(1)}.
Here is the plan of the proof. First, we estimate how close the $K$%
-functions of $\mu^{\left( n\right) }$ and $\Phi$ are to each other. Then we
note that the Cauchy transforms of $\mu^{\left( n\right) }$ and $\Phi$ can
be found from their functional equations:%
\begin{equation*}
K_{n}\left( G_{n}\left( z\right) \right) =z, 
\end{equation*}
and 
\begin{equation*}
K_{\Phi}\left( G_{\Phi}\left( z\right) \right) =z, 
\end{equation*}
where $K_{n}$ and $K_{\Phi}$ denote the $K$-functions of $\mu^{\left(
n\right) }$ and $\Phi.$ From the previous step we know that $K_{n}\left(
z\right) $ is close to $K_{\Phi}\left( z\right) .$ Our goal is to show that
this implies that $G_{n}\left( z\right) $ is close to $G_{\Phi}\left(
z\right) .$

If we introduce an extra parameter, $t,$ then we can include these
functional equations in a parametric family:%
\begin{equation}
K_{t}\left( G_{t}\left( z\right) \right) =z.  \label{functional_equation}
\end{equation}
Parameter $t=0$ corresponds to $\Phi$ and $t=1$ to $\mu^{\left( n\right) }$.
Next, we fix $z$ and consider $G_{t}$ as a function of $t.$ We develop this
function in a power series in $t$:%
\begin{equation*}
G_{t}=G_{\Phi}+\sum_{k=1}^{\infty}c_{k}t^{k}, 
\end{equation*}
where $c_{k}$ are functions of $z.$ Then we estimate $I_{k}$ for each $%
k\geq1,$ where 
\begin{equation*}
I_{k}=\int\limits_{-8}^{8}\left\vert c_{k}\left( u+iv\right) \right\vert du. 
\end{equation*}
Then,%
\begin{equation*}
\int\limits_{-8}^{8}\left\vert G_{n}\left( u+iv\right) -G_{\Phi}\left(
u+iv\right) \right\vert du\leq\sum_{k=1}^{\infty}I_{k}, 
\end{equation*}
and our estimates of $I_{k}$ allow us to prove the claim of Lemma \ref%
{Cauchy_tr_closeness copy(1)}.

\subsection{Speed of convergence of K-functions}

In this section, we start with some preliminary facts about functional
inverses of holomorphic functions. Then we derive an estimate for the speed
of convergence of the $K$-functions of $\mu^{\left( n\right) }$ and the
semicircle law.

First, recall the Lagrange formula for the functional inversion. By a
function holomorphic in a domain $D$ we mean a function which is bounded and
complex-differentiable in $D$.

\begin{lemma}[Lagrange's inversion formula]
\label{Lagranges_series_around0} Suppose $f$ is a function of a complex
variable, which is holomorphic in a neighborhood of $z_{0}$ and suppose that 
$f\left( z_{0}\right) =0$ and $f^{\prime }(z_{0})\neq 0.$ Then the
functional inverse of $f\left( z\right) $ is well defined in a neighborhood
of $0$ and the Taylor series of the inverse is given by the following
formula:%
\begin{equation*}
f^{-1}\left( u\right) =z_{0}+\sum_{k=1}^{\infty }\left[ \frac{1}{k}\mathrm{res}_{z=z_{0}}\frac{1}{f(z)^{k}}\right] u^{k},
\end{equation*}%
where $\mathrm{res}_{z=z_{0}}$ denotes the Cauchy residual at point $z_{0}$.
\end{lemma}

For proof see Theorems II.3.2 and II.3.3 in \shortciteN{markushevich77}, or
Section 7.32 in \shortciteN{whittaker_watson27}. We also need the following
modification of the Lagrange formula.

\begin{lemma}
\label{Lagranges_series} Suppose $G$ is a function of a complex variable,
which is holomorphic in a neighborhood of $z_{0}=\infty$ and has the
expansion 
\begin{equation*}
G(z)=\frac{1}{z}+\frac{a_{1}}{z^{2}}+..., 
\end{equation*}
converging for all sufficiently large $z.$ Define $g(z)=G(1/z).$ Then the
functional inverse of $G\left( z\right) $ is well defined in a neighborhood
of $0$ and the Laurent series of the inverse is given by the following
formula:%
\begin{equation*}
G^{-1}\left( w\right) =\frac{1}{w}+a_{1}-\sum_{n=1}^{\infty}\left[ \frac {1}{%
n}\frac{1}{2\pi i}\oint_{\partial\gamma}\frac{dz}{z^{2}g(z)^{n}}\right]
w^{n}, 
\end{equation*}
where $\gamma$ is a closed disc around $0$ in which $g(z)$ has only one zero$%
.$
\end{lemma}

\textbf{Proof: }Let $\gamma$ be a disc around $0$ in which $g(z)$ has only
one zero. This disk exists because $g\left( 0\right) =0,$ and $g\left(
z\right) $ is analytical in a neighborhood of $0$ and has a non-zero
derivative at $0.$ Let 
\begin{equation*}
r_{w}=\frac{1}{2}\inf_{z\in\partial\gamma}\left\vert g\left( z\right)
\right\vert . 
\end{equation*}
Then $r_{w}>0$ by our assumption about $\gamma.$ We can apply Rouch\'{e}'s
theorem and conclude that the equation $g\left( z\right) -w=0$ has only one
solution inside $\gamma$ if $\left\vert w\right\vert \leq r_{w}.$ Let us fix
such $w$ that $\left\vert w\right\vert \leq r_{w}$. Inside $\gamma$, the
function 
\begin{equation*}
\frac{g^{\prime}(z)}{z\left( g(z)-w\right) } 
\end{equation*}
has a pole at $z=1/G^{-1}(w)$ with the residual $G^{-1}(w)$ and a pole at $%
z=0$ with the residual $-1/w$. Consequently, we can write:%
\begin{equation*}
G^{-1}\left( w\right) =\frac{1}{2\pi i}\oint_{\partial\gamma}\frac {%
g^{\prime}(z)dz}{z\left( g(z)-w\right) }+\frac{1}{w}. 
\end{equation*}
The integral can be re-written as follows:%
\begin{align*}
\oint_{\partial\gamma}\frac{g^{\prime}(z)dz}{z\left( g(z)-w\right) } &
=\oint_{\partial\gamma}\frac{g^{\prime}(z)}{zg(z)}\frac{1}{1-\frac{w}{%
g\left( z\right) }}dz \\
& =\sum_{n=0}^{\infty}\oint_{\partial\gamma}\frac{g^{\prime}(z)dz}{%
zg(z)^{n+1}}w^{n}.
\end{align*}
For $n=0$ we calculate 
\begin{equation*}
\frac{1}{2\pi i}\oint_{\partial\gamma}\frac{g^{\prime}(z)dz}{zg(z)}=a_{1}, 
\end{equation*}
Indeed, the only pole of the integrand is at $z=0$ and it has order two. The
corresponding residual can be computed from the series expansion for $g(z)$:%
\begin{align*}
\mathrm{res}_{z=0}\frac{g^{\prime}(z)dz}{zg(z)} & =\left. \frac{d}{dz}\frac{%
z^{2}\left( 1+2a_{1}z+...\right) }{z\left( z+a_{1}z^{2}+..\right) }%
\right\vert _{z=0} \\
& =\left. \frac{d}{dz}\frac{1+2a_{1}z+...}{1+a_{1}z+..}\right\vert
_{z=0}=a_{1}.
\end{align*}
For $n>0$ we integrate by parts:%
\begin{equation*}
\frac{1}{2\pi i}\oint_{\partial\gamma}\frac{g^{\prime}(z)dz}{zg(z)^{n+1}}=-%
\frac{1}{2\pi i}\frac{1}{n}\oint_{\partial\gamma}\frac{dz}{z^{2}g(z)^{n}}. 
\end{equation*}
QED.

Our first application of the Lagrange formula is the following estimate on
the convergence radius of power series for $K$-functions.

\begin{lemma}
Suppose that the measure $\mu$ is supported on interval $\left[ -L,L\right]
, $ and $K\left( z\right) $ denotes the functional inverse of its Cauchy
transform $G\left( z\right) .$ Then the Laurent series of $K\left( z\right) $
converge in the area $\Omega=\left\{ z:0<\left\vert z\right\vert <\left(
4L\right) ^{-1}\right\} .$
\end{lemma}

\textbf{Proof:} Let us apply Lemma \ref{Lagranges_series} to $G\left(
z\right) $ with circle $\gamma$ having radius $\left( 2L\right) ^{-1}.$ We
need to check that $g(z)=:G\left( 1/z\right) $ has only one zero inside this
circle. This holds because 
\begin{equation*}
g(z)=z\left( 1+a_{2}z^{2}+a_{3}z^{3}+...\right) , 
\end{equation*}
and inside $\left\vert z\right\vert \leq\left( 2L\right) ^{-1}$ we can
estimate: 
\begin{equation}
\left\vert a_{2}z^{2}+a_{3}z^{3}+...\right\vert \leq L^{2}\left( \frac{1}{2L}%
\right) ^{2}+L^{3}\left( \frac{1}{2L}\right) ^{3}+...=\frac{1}{2},
\label{K-convergence}
\end{equation}
and an application of Rouch\'{e}'s theorem shows that $g\left( z\right) $
has only one zero inside this circle.

Another consequence of the estimate (\ref{K-convergence}) is that on the
circle $\left\vert z\right\vert =\left( 2L\right) ^{-1}$ 
\begin{equation*}
\left\vert g\left( z\right) \right\vert \geq\left\vert z\right\vert
/2=1/\left( 4L\right) . 
\end{equation*}

By Lemma \ref{Lagranges_series} the coefficients in the series for the
inverse of $G\left( z\right) $ are 
\begin{equation*}
b_{k}=\frac{1}{2\pi ik}\oint_{\partial\gamma}\frac{dz}{z^{2}g(z)^{k}}, 
\end{equation*}
and we can estimate them as%
\begin{equation}
\left\vert b_{k}\right\vert \leq\frac{2L}{k}\left( 4L\right) ^{k}.
\label{b_estimates}
\end{equation}
This implies that the radius of convergence of power series for $K\left(
z\right) $ is at least $\left( 4L\right) ^{-1}.$ QED.

Let $K_{n}\left( z\right) $ denote the $K$-function of $\mu^{\left( n\right)
}$. For the semicircle law the $K$-function is $K_{\Phi}\left( z\right)
=z^{-1}+z$. Define $\varphi_{n}\left( z\right) =K_{n}\left( z\right)
-z-z^{-1}.$

\begin{lemma}
\label{phi_size} Suppose $\mu$ has zero mean and unit variance and is
supported on $\left[ -L,L\right] $. Then the function $\varphi_{n}\left(
z\right) $ is holomorphic in $\left\vert z\right\vert \leq\sqrt{n}/\left(
8L\right) $ and 
\begin{equation*}
\left\vert \varphi_{n}\left( z\right) \right\vert \leq32L^{3}\frac {%
\left\vert z\right\vert ^{2}}{\sqrt{n}}. 
\end{equation*}
\end{lemma}

\textbf{Proof:} The measure $\mu^{\left( n\right) }$ is the $n$-time
convolution of the measure $\widetilde{\mu}\left( dx\right) =:\mu\left( 
\sqrt{n}dx\right) $ with itself. Therefore, $K_{n}\left( z\right) =nK_{%
\widetilde{\mu}}\left( z\right) -\left( n-1\right) z^{-1}$. Since $%
\widetilde{\mu}$ is supported on $\left[ -L/\sqrt{n},L/\sqrt{n}\right] $, we
can estimate $K_{\widetilde{\mu}}\left( z\right) -\frac{1}{z}-\frac{1}{n}z$
inside the circle $\left\vert z\right\vert =\sqrt{n}/\left( 8L\right) $
using the estimates for coefficients of $K_{\widetilde{\mu}}\left( z\right) $
in (\ref{b_estimates}) and changing $L$ to $L/\sqrt{n}$ in these estimates:%
\begin{align*}
\left\vert K_{\widetilde{\mu}}\left( z\right) -\frac{1}{z}-\frac{1}{n}%
z\right\vert & =\sum_{k=2}^{\infty}b_{k}z^{k}\leq\frac{2L}{\sqrt{n}}%
\sum_{k=2}^{\infty}\frac{1}{k}\left( \frac{4L}{\sqrt{n}}\right)
^{k}\left\vert z\right\vert ^{k} \\
& \leq32\left( \frac{L}{\sqrt{n}}\right) ^{3}\left\vert z\right\vert
^{2}\sum_{k=2}^{\infty}\frac{1}{k2^{k-1}} \\
& \leq32\left( \frac{L}{\sqrt{n}}\right) ^{3}\left\vert z\right\vert ^{2}.
\end{align*}
Note that we used the assumption about the mean and variance of the measure $%
\mu$ in the first line by setting $b_{1}=0$ and $b_{1}=1/n.$

Using the summation formula (\ref{summation_property}) for $K$-functions, we
further obtain: 
\begin{equation*}
\left\vert K_{n}\left( z\right) -\frac{1}{z}-z\right\vert \leq32\frac{L^{3}}{%
\sqrt{n}}\left\vert z\right\vert ^{2}. 
\end{equation*}
QED.

Lemma \ref{phi_size} shows that as $n$ grows, the radius of the convergence
area of $\varphi_{n}\left( z\right) ,$ and therefore of $K_{n}\left(
z\right) ,$ grows proportionally to $\sqrt{n}.$ In particular, the radius of
convergence will eventually cover every bounded domain. Lemma \ref{phi_size}
also establishes the rate of convergence of $K_{n}\left( z\right) $ to its
limit $K_{\Phi}\left( z\right) =z^{-1}+z.$

\subsection{Useful Estimates}

Suppose $G_{\Phi}\left( z\right) $ is the Cauchy transform of the semicircle
distribution.

\begin{lemma}
\label{GSC_size} 1) $\left\vert G_{\Phi}\left( z\right) \right\vert \leq1$
if $\mathrm{Im}z>0;$\newline
2) $\left\vert z-2G_{\Phi}\left( z\right) \right\vert \geq2\sqrt{\mathrm{Im}z}$
if $\mathrm{Im}z\in\left( 0,2\right) .$
\end{lemma}

\textbf{Proof:} $G_{\Phi}\left( z\right) =\left( z-\sqrt{z^{2}-4}\right) /2. 
$ If $z=u+iv$ and $v$ is fixed, then the maximum of $\left\vert G_{\Phi
}\left( z\right) \right\vert $ is reached for $u=0.$ Then $\left\vert
G_{\Phi}\left( iv\right) \right\vert =\left( \sqrt{v^{2}+4}-v\right) /2$ and 
$\sup\left\vert G_{\Phi}\left( iv\right) \right\vert =1.$

Next, $\left\vert z-2G_{\Phi}\left( z\right) \right\vert =\left\vert \sqrt{%
u^{2}-v^{2}-4+i2uv}\right\vert .$ If $v$ is in $\left( 0,2\right) $ and
fixed, the minimum of this expression is reached for $u=\pm\sqrt{4-v^{2}}$
and equals $2\sqrt{v}.$ QED.

\begin{lemma}
\label{Phi_GSC_size} If $n\geq64L^{2}$ and $\mathrm{Im}z>0,$ then we have: 
\begin{equation*}
\left\vert \varphi_{n}\left( G_{\Phi}\left( z\right) \right) \right\vert \leq%
\frac{32L^{3}}{\sqrt{n}}. 
\end{equation*}
\end{lemma}

\textbf{Proof:} This Lemma follows directly from Lemmas \ref{phi_size} and %
\ref{GSC_size}. QED.

\subsection{Functional equation for the Cauchy transform}

Let $G_{n}\left( z\right) $ denote the Cauchy transform of $\mu^{\left(
n\right) }.$ Let us write the following functional equation:%
\begin{equation}
G\left( t,z\right) +\frac{1}{G\left( t,z\right) }+t\varphi_{n}\left( G\left(
t,z\right) \right) =z,  \label{inversion_equation}
\end{equation}
where $t$ is a complex parameter. For $t=0$ the solution is $G_{\Phi}\left(
z\right) $, and for $t=1$ the solution is $G_{n}\left( z\right) .$ Assume
that $\varphi_{n}\left( z\right) $ is not identically zero. (If it is, then $%
\mu^{\left( n\right) }$ is semicircle and $d\left( \mu^{\left( n\right)
},\nu\right) =0.$) Let us write equation (\ref{inversion_equation}) as 
\begin{equation}
t=\frac{zG-G^{2}-1}{G\varphi_{n}\left( G\right) }.
\label{inversion_equation2}
\end{equation}
We can think about $z$ as a fixed complex parameter and about $t$ as a
function of the complex variable $G$, i.e., $t=f\left( G\right) .$ Suppose $%
\varphi_{n}\left( G_{\Phi}\left( z\right) \right) $ does not equal zero for
a given value of $z$. (This holds for all but a countable number of values
of parameter $z$.) Then, as a function of $G$, $f$ is holomorphic in a
neighborhood of $G_{\Phi}\left( z\right) .$\ What we would like to do is to
invert this function $f$ and write $G=f^{-1}\left( t\right) .$ In particular
we would like to develop $f^{-1}\left( t\right) $ in a series of $t$ around $%
t=0.$ Then we would be able to estimate $\left\vert f^{-1}\left( 1\right)
-f^{-1}\left( 0\right) \right\vert $, which is equal to $\left\vert
G_{n}\left( z\right) -G_{\Phi}\left( z\right) \right\vert .$ To perform this
inversion, we use the Lagrange formula in Lemma \ref%
{Lagranges_series_around0}.

Assume that $z$ is fixed, and let us write $G$ instead of $G\left( z\right) $
and $G_{\Phi}$ instead of $G_{\Phi}\left( z\right) .$ By Lemma \ref%
{Lagranges_series_around0}, we can write the solution of (\ref%
{inversion_equation2}) as 
\begin{equation}
G=G_{\Phi}+\sum_{k=1}^{\infty}c_{k}t^{k},  \label{formula_G_expansion}
\end{equation}
where 
\begin{equation}
c_{k}=\frac{1}{k}\mathrm{res}_{G=G_{\Phi}}\left( \frac{G\varphi _{n}\left(
G\right) }{zG-G^{2}-1}\right) ^{k}.  \label{formula_ck_integral_expression}
\end{equation}

We aim to estimate $I_{k}=:$ $\int_{-8}^{8}\left\vert c_{k}\left(
u+iv\right) \right\vert du.$ In particular, we will show that for any $%
v\in\left( 0,1\right) ,$ $I_{1}=O\left( n^{-1/2}\right) .$ In addition, we
will show that if $v=b/\sqrt{n}$ for a suitably chosen $b,$ then $%
\sum_{k=2}^{\infty}I_{k}=o\left( n^{-1/2}\right) .$ This information is
sufficient for a good estimate of $\int_{-8}^{8}\left\vert G_{n}\left(
u+iv\right) -G_{\Phi}\left( u+iv\right) \right\vert du.$

Let us consider first the case of $k=1.$ Then 
\begin{align*}
c_{1} & =\frac{G_{\Phi}\varphi_{n}\left( G_{\Phi}\right) }{G_{2}-G_{\Phi}} \\
& =\frac{G_{\Phi}\varphi_{n}\left( G_{\Phi}\right) }{\sqrt{z^{2}-4}},
\end{align*}
where $G_{2}$ denotes the root of the equation $G^{2}-zG+1=0,$ which is
different from $G_{\Phi}.$ Therefore, if $z=u+iv,$ then we can calculate: 
\begin{align*}
\left\vert c_{1}\right\vert & =\frac{\left\vert G_{\Phi}\right\vert
\left\vert \varphi_{n}\left( G_{\Phi}\right) \right\vert }{\left[ \left(
u^{2}-4\right) ^{2}+2\left( u^{2}+4\right) v^{2}+v^{4}\right] ^{1/4}} \\
& \leq\frac{32L^{3}}{\sqrt{n}}\frac{1}{\left[ \left( u^{2}-4\right)
^{2}+2\left( u^{2}+4\right) v^{2}+v^{4}\right] ^{1/4}},
\end{align*}
where the last inequality holds by Lemma \ref{Phi_GSC_size} for all $%
n\geq64L^{2}.$

\begin{lemma}
For every $v\in\left( 0,1\right) ,$%
\begin{equation*}
\int \limits_{-8}^{8}\left[ \left( u^{2}-4\right) ^{2}+2\left(
u^{2}+4\right) v^{2}+v^{4}\right] ^{-1/4}du<24. 
\end{equation*}
\end{lemma}

\textbf{Proof:} Let us make substitution $x=u^{2}-4.$ Then we get:%
\begin{align*}
J & =\int \limits_{-4}^{60}\left[ x^{2}+2\left( x+8\right) v^{2}+v^{4}%
\right] ^{-1/4}\frac{dx}{\sqrt{x+4}} \\
& \leq\int \limits_{-4}^{60}\frac{1}{\left( x^{2}+v^{2}\right) ^{1/4}}\frac{%
dx}{\sqrt{x+4}}.
\end{align*}
Now we divide the interval of integration in two parts and write: 
\begin{align*}
J & \leq\int \limits_{-4}^{-2}...+\int\limits_{-2}^{60}... \\
& \leq\frac{1}{\sqrt{2}}\int \limits_{-4}^{-2}\frac{dx}{\sqrt{x+4}}+\frac{2%
}{\sqrt{2}}\int\limits_{0}^{60}\frac{dx}{x^{1/2}} \\
& =2+\frac{4}{\sqrt{2}}\sqrt{60}<24.
\end{align*}
QED.

\begin{corollary}
\label{Corollary_I1}For every $v\in\left( 0,1\right) $ and all $n\geq
64L^{2},$ it is true that 
\begin{equation*}
I_{1}=:\int\limits_{-8}^{8}\left\vert c_{1}\left( u+iv\right) \right\vert
du\leq\frac{768L^{3}}{\sqrt{n}} 
\end{equation*}
\end{corollary}

Now we estimate $c_{k}$ in (\ref{formula_G_expansion}) for $k\geq2$. Define
function $f_{k}\left( G\right) $ by the formula%
\begin{equation*}
f_{k}\left( G\right) =:\left( \frac{G\varphi_{n}\left( G\right) }{G_{2}-G}%
\right) ^{k}, 
\end{equation*}
where $G_{2}$ denotes the root of the equation $G^{2}-zG+1=0,$ which is
different from $G_{\Phi}.$ Then formula (\ref{formula_ck_integral_expression}%
) implies that $kc_{k}$ equal to the coefficient before $\left( G-G_{\Phi
}\right) ^{k-1}$ in the expansion of $f_{k}\left( G\right) $ in power series
of $\left( G-G_{\Phi}\right) $. To estimate this coefficient, we will use
the Cauchy inequality:%
\begin{equation*}
\left\vert kc_{k}\right\vert \leq\frac{M_{k}\left( r\right) }{r^{k-1}}, 
\end{equation*}
where $M_{k}\left( r\right) $ is the maximum of $\left\vert f_{k}\left(
G\right) \right\vert $ on the circle $\left\vert G-G_{\Phi}\right\vert =r.$

We will use $r=\sqrt{v}$ and our first goal is to estimate $M_{k}\left( 
\sqrt{v}\right) .$

\begin{lemma}
Let $z=u+iv$ and suppose that $v\in\left( 0,1\right) .$ If $n\geq256L^{2},$
then 
\begin{equation*}
M_{k}\left( v\right) \leq\left[ \frac{512L^{3}}{\sqrt{n}}\frac{1}{\left(
\left( u^{2}-4\right) ^{2}+2u^{2}v^{2}\right) ^{1/4}}\right] ^{k}. 
\end{equation*}
\end{lemma}

\textbf{Proof:} Note that $\left\vert G\right\vert \leq\left\vert G_{\Phi
}\right\vert +\sqrt{v}$ and therefore $\left\vert G\right\vert \leq2$
provided that $v\in\left( 0,1\right) .$ Then Lemma \ref{phi_size} implies
that if $n\geq256L^{2},$ then $\varphi_{n}\left( G\right) $ is well defined
and $\left\vert \varphi_{n}\left( G\right) \right\vert \leq128L^{3}/\sqrt{n}%
. $ It remains to estimate $\left\vert G_{2}-G\right\vert $ from below. If
we write $G=G_{\Phi}+e^{i\theta}\sqrt{v},$ then we have 
\begin{align*}
\left\vert G_{2}-G\right\vert & =\left\vert \sqrt{z^{2}-4}-e^{i\theta}\sqrt{v%
}\right\vert \\
& \geq\left\vert \sqrt{z^{2}-4}\right\vert -\sqrt{v} \\
& =\left( \left( u^{2}-4\right) ^{2}+2\left( u^{2}+4\right)
v^{2}+v^{4}\right) ^{1/4}-\sqrt{v} \\
& >\left( \left( u^{2}-4\right) ^{2}+2\left( u^{2}+4\right) v^{2}\right)
^{1/4}-\sqrt{v}>0.
\end{align*}

From the concavity of function $t^{1/4}$ it follows that for positive $A$
and $B$ the following inequality holds: 
\begin{equation*}
\left[ 8\left( A+B\right) \right] ^{1/4}-A^{1/4}\geq B^{1/4}. 
\end{equation*}
Using $v^{2}$ as $A,$ and $\left[ \left( u^{2}-4\right) ^{2}+2u^{2}v^{2}%
\right] /8$ as $B,$ we can write this inequality as follows:%
\begin{equation*}
\left( \left( u^{2}-4\right) ^{2}+2\left( u^{2}+4\right) v^{2}\right) ^{1/4}-%
\sqrt{v}\geq\frac{1}{8^{1/4}}\left( \left( u^{2}-4\right)
^{2}+2u^{2}v^{2}\right) ^{1/4}>0. 
\end{equation*}

Therefore 
\begin{equation*}
M_{k}\left( v\right) \leq\left[ \frac{512L^{3}}{\sqrt{n}}\frac{1}{\left(
\left( u^{2}-4\right) ^{2}+2u^{2}v^{2}\right) ^{1/4}}\right] ^{k}. 
\end{equation*}
QED.

\begin{corollary}
\label{corollary_ck_estimate}For every $v\in\left( 0,1\right) ,$ $k\geq2,$
and all $n\geq256L^{2},$ it is true that 
\begin{equation*}
\left\vert kc_{k}\left( u+iv\right) \right\vert \leq v^{-\frac{k-1}{2}}\left[
\frac{512L^{3}}{\sqrt{n}}\frac{1}{\left( \left( u^{2}-4\right)
^{2}+2u^{2}v^{2}\right) ^{1/4}}\right] ^{k} 
\end{equation*}
\end{corollary}

Now we want to estimate integrals of $\left\vert c_{k}\left( u+iv\right)
\right\vert $ when $u$ changes from $-8$ to $8.$ The cases of $k=2$ and $k>2$
are slightly different and we treat them separately.

Let 
\begin{equation*}
I_{k}=:\int\limits_{-8}^{8}\left\vert c_{k}\left( u+iv\right) \right\vert du 
\end{equation*}

\begin{lemma}
If $v\in\left( 0,1\right) $ and $n\geq256L^{2}$, then i) 
\begin{equation*}
I_{2}\leq\frac{\log\left( 60/v\right) }{\sqrt{v}}\frac{2^{19}L^{6}}{n}, 
\end{equation*}
and ii) if $k>2,$ then%
\begin{equation*}
I_{k}\leq\frac{12}{k}v^{3/2}\left( \frac{512L^{3}}{v\sqrt{n}}\right) ^{k}. 
\end{equation*}
\end{lemma}

\textbf{Proof: }Using Corollary \ref{corollary_ck_estimate}, we write:%
\begin{equation*}
I_{k}=:\int\limits_{8}^{8}\left\vert c_{k}\left( u+iv\right) \right\vert
du\leq\frac{1}{k}\frac{1}{v^{\left( k-1\right) /2}}\left( \frac{512L^{3}}{%
\sqrt{n}}\right) ^{k}\int\limits_{-8}^{8}\frac{du}{\left( \left(
u^{2}-4\right) ^{2}+2u^{2}v^{2}\right) ^{k/4}}. 
\end{equation*}
After substitution $x=u^{2}-4,$ the integral in the right-hand side of the
inequality can be re-written as 
\begin{equation*}
J_{k}=:\int\limits_{-4}^{60}\frac{1}{\left( x^{2}+2\left( x+4\right)
v^{2}\right) ^{k/4}}\frac{dx}{\sqrt{x+4}}. 
\end{equation*}
We divide the interval of integration into two portions and write the
following inequality:%
\begin{align*}
J_{k} & \leq\int\limits_{-4}^{-1}\frac{dx}{\sqrt{x+4}}+\int\limits_{-1}^{60}%
\frac{dx}{\left( x^{2}+2\left( x+4\right) v^{2}\right) ^{k/4}} \\
& \leq2\sqrt{3}+\int\limits_{-1}^{60}\frac{dx}{\left( x^{2}+v^{2}\right)
^{k/4}}.
\end{align*}
If we use substitution $s=x/v,$ then we can write:%
\begin{align*}
\int\limits_{-1}^{60}\frac{dx}{\left( x^{2}+v^{2}\right) ^{k/4}} &
=\int\limits_{-1/v}^{60/v}\frac{ds}{v^{\frac{k}{2}-1}\left( 1+s^{2}\right)
^{k/4}} \\
& \leq\frac{2}{v^{\frac{k}{2}-1}}\int\limits_{0}^{60/v}\frac{ds}{\left(
1+s^{2}\right) ^{k/4}}.
\end{align*}
We again separate the interval of integration in two parts and write:%
\begin{equation*}
\frac{2}{v^{\frac{k}{2}-1}}\int\limits_{0}^{60/v}\frac{ds}{\left(
1+s^{2}\right) ^{k/4}}\leq\frac{2}{v^{\frac{k}{2}-1}}\left[ \int
\limits_{0}^{1}ds+\int\limits_{1}^{60/v}\frac{ds}{s^{k/2}}\right] . 
\end{equation*}
Here we have two different cases. If $k=2,$ then we evaluate the integrals
as $1+\log\left( 60/v\right) .$ Therefore, 
\begin{align*}
J_{2} & \leq2\sqrt{3}+2+2\log\left( 60/v\right) \\
& \leq4\log\left( 60/v\right) .
\end{align*}
Hence, 
\begin{equation*}
I_{2}\leq\frac{\log\left( 60/v\right) }{\sqrt{v}}\frac{2^{19}L^{6}}{n}. 
\end{equation*}

If $k>2,$ then we have:%
\begin{align*}
\frac{2}{v^{\frac{k}{2}-1}}\left[ \int\limits_{0}^{1}ds+\int%
\limits_{1}^{60/v}\frac{ds}{s^{k/2}}\right] & =\frac{2}{v^{\frac{k}{2}-1}}%
\left[ 1+\frac{1}{-\frac{k}{2}+1}\left( \left( \frac{60}{v}\right) ^{-\frac{k%
}{2}+1}-1\right) \right] \\
& \leq\frac{2}{v^{\frac{k}{2}-1}}\left( 1+\frac{1}{\frac{k}{2}-1}\right) . \\
& \leq\frac{6}{v^{\frac{k}{2}-1}}.
\end{align*}
Therefore, 
\begin{equation*}
J_{k}\leq2\sqrt{3}+\frac{6}{v^{\frac{k}{2}-1}}\leq\frac{12}{v^{\frac{k}{2}-1}%
}, 
\end{equation*}
and 
\begin{align*}
I_{k} & \leq\frac{1}{k}\frac{1}{v^{\left( k-1\right) /2}}\left( \frac{%
512L^{3}}{\sqrt{n}}\right) ^{k}\frac{12}{v^{\frac{k}{2}-1}} \\
& \leq\frac{12}{k}v^{3/2}\left( \frac{512L^{3}}{v\sqrt{n}}\right) ^{k}.
\end{align*}

QED.

\begin{corollary}
\label{Corollary_I2}If $v=1024L^{3}/\sqrt{n},$ and $n\geq256L^{2}$ then 
\begin{equation*}
I_{2}\leq2^{14}L^{9/2}\frac{\log\left( \frac{15}{256L^{3}}\sqrt{n}\right) }{%
n^{3/4}}. 
\end{equation*}
In particular, $I_{2}=o\left( n^{-1/2}\right) $ as $n\rightarrow\infty.$
\end{corollary}

Now we address the case when $k>2.$

\begin{corollary}
\label{Corollary_I3}Suppose $v=1024L^{3}/\sqrt{n},$ and $n\geq256L^{2}$.
Then 
\begin{equation*}
\sum\limits_{k=3}^{\infty}I_{k}\leq\frac{3}{2}v^{3/2}=1536L^{9/2}\frac {1}{%
n^{3/4}}. 
\end{equation*}
In particular, $\sum\nolimits_{k=2}^{\infty}I_{k}=o\left( n^{-1/2}\right) $
as $n\rightarrow\infty.$
\end{corollary}

Joining results of Corollaries \ref{Corollary_I1}, \ref{Corollary_I2}, and %
\ref{Corollary_I3}, we get the following result.

\begin{lemma}
\label{Cauchy_tr_closeness}Suppose $v=1024L^{3}/\sqrt{n}.$ Then for all
sufficiently large $n,$ we have the following estimate:%
\begin{equation*}
\int\limits_{-8}^{8}\left\vert G_{n}\left( u+iv\right) -G_{\Phi}\left(
u+iv\right) \right\vert du\leq\frac{2^{10}L^{3}}{\sqrt{n}}. 
\end{equation*}
\end{lemma}

\textbf{Proof:} From formula (\ref{formula_G_expansion}) we have 
\begin{equation*}
\left\vert G_{n}-G_{\Phi}\right\vert \leq\sum_{k=1}^{\infty}\left\vert
c_{k}\right\vert . 
\end{equation*}
Since the series has only positive terms, we can integrate it term by term
and write:%
\begin{align*}
\int\limits_{-8}^{8}\left\vert G_{n}\left( u+iv\right) -G_{\Phi}\left(
u+iv\right) \right\vert du &
\leq\sum_{k=1}^{\infty}\int\limits_{-8}^{8}\left\vert c_{k}\left(
u+iv\right) \right\vert du \\
& \leq\frac{768L^{3}}{\sqrt{n}}+o\left( n^{-1/2}\right) \\
& \leq\frac{2^{10}L^{3}}{\sqrt{n}}
\end{align*}
for all sufficiently large $n.$ QED.

Lemma \ref{Cauchy_tr_closeness} is identical to Lemma \ref%
{Cauchy_tr_closeness copy(1)} and its proof completes the proof of Theorem %
\ref{theorem_main}.

\section{Example}

Consider a binomial measure: $\mu\left\{ -1/p\right\} =p$ and $\mu\left\{
1/q\right\} =q\equiv1-p.$ This is a zero-mean measure with the variance
equal to $\left( pq\right) ^{-1}.$ Let $\mu^{\left( n\right) }\left(
dx\right) =\mu\boxplus...\boxplus\mu\left( \sqrt{\frac{n}{pq}}dx\right) $
and let $\mathcal{F}_{n}\left( x\right) $ be the distribution function
corresponding to $\mu^{\left( n\right) }.$

\begin{proposition}
If $p\neq q,$ then there exist such positive constants $C_{1}$ $\ $\ and $%
C_{2}$ that 
\begin{equation*}
C_{1}n^{-1/2}\leq\sup_{x}\left\vert \mathcal{F}_{n}\left( x\right)
-\Phi\left( x\right) \right\vert \leq C_{2}n^{-1/2} 
\end{equation*}
for every $n.$
\end{proposition}

\textbf{Proof:} From the Voiculescu addition formula and the Stieltjes
inversion formula, it is easy to compute the density of the distribution of $%
\mu^{\left( n\right) }$:%
\begin{equation*}
f_{n}\left( x\right) =\frac{1}{2\pi}\frac{\sqrt{4-x^{2}+2\frac{p-q}{\sqrt{pqn%
}}x-\frac{1}{pqn}}}{\left( 1+\frac{x}{\sqrt{nq/p}}\right) \left( 1-\frac{x}{%
\sqrt{np/q}}\right) }, 
\end{equation*}
if the square root is real, and if not, $f_{n}\left( x\right) =0.$ We
compare this distribution with the semicircle distribution, which has the
following density:%
\begin{equation*}
\phi\left( x\right) =\frac{1}{2\pi}\sqrt{4-x^{2}}\chi_{\lbrack-2,2]}\left(
x\right) . 
\end{equation*}
More precisely, we seek to estimate 
\begin{equation*}
\sup_{x}\left\vert \int_{-\infty}^{x}\left( f_{n}\left( t\right) -\phi\left(
t\right) \right) dt\right\vert . 
\end{equation*}

The support of $f_{n}$ is $\left[ -2\sqrt{1-n^{-1}}+cn^{-1/2},\text{ }2\sqrt{%
1-n^{-1}}+cn^{-1/2}\right] ,$ where $c=\left( p-q\right) /\sqrt {pq}.$
Suppose in the following that $p>q$ and introduce the new variable $u=x+2%
\sqrt{1-n^{-1}}-c/\sqrt{n}.$ Then, 
\begin{align*}
2\pi f_{n}\left( u\right) & =\sqrt{4u\sqrt{1-n^{-1}}-u^{2}} \\
& \sim\sqrt{4u-u^{2}},
\end{align*}
where the asymptotic equivalence is for $u$ fixed and $n\rightarrow\infty$
and we omit all terms that are $o\left( n^{-1/2}\right) .$ Similarly, 
\begin{align*}
2\pi\phi\left( x\right) & =\sqrt{4-\left( u-2\left( 1-n^{-1}\right)
+cn^{-1/2}\right) ^{2}} \\
& \sim\sqrt{4u-u^{2}+4cn^{-1/2}-2cun^{-1/2}} \\
& =\sqrt{4u-u^{2}}\sqrt{1+2c\frac{1}{\sqrt{n}u}\frac{2-u}{4-u}}.
\end{align*}
Consequently, 
\begin{align*}
\phi\left( u\right) -f_{n}\left( u\right) & \sim\sqrt{4u-u^{2}}\left[ \sqrt{%
1+2c\frac{1}{\sqrt{n}u}\frac{2-u}{4-u}}-1\right] \\
& \sim\sqrt{4u-u^{2}}c\frac{1}{\sqrt{n}u}\frac{2-u}{4-u} \\
& =c\frac{1}{\sqrt{nu}}\frac{2-u}{\sqrt{4-u}}.
\end{align*}
After integrating we get:%
\begin{equation*}
\left\vert \int_{0}^{x}\left( f_{n}\left( u\right) -\phi\left( u\right)
\right) du\right\vert \sim c\frac{1}{\sqrt{n}}f\left( x\right) , 
\end{equation*}
where $f\left( x\right) $ is a continuous positive bounded function. From
this expression it is clear that $\sup_{x}\left\vert
\int_{-\infty}^{x}\left( f_{n}\left( t\right) -\phi\left( t\right) \right)
dt\right\vert $ has the order of $n^{-1/2}$ provided that $p\neq q.$ QED.

This example shows that the rate of $n^{-1/2}$ in Theorem \ref{theorem_main}
cannot be improved without further restrictions on measures. It would be
interesting to extend Theorem \ref{theorem_main} to measures with unbounded
support or relate the constant in the inequality to moments of the convolved
measures, similar to the classical case.

\bibliographystyle{CHICAGO}
\bibliography{comtest}

\end{document}